\documentclass[12pt]{article}
\usepackage[utf8]{inputenc}

\usepackage{newtxtext}

\usepackage{amsfonts}

\usepackage{mdframed}

\usepackage{amssymb, amsmath}

\usepackage[round]{natbib}

\usepackage{setspace} 

\usepackage{footmisc}

\newtheorem{definition}{Definition}

\author{Claudio Ternullo, Isabella Fascitiello}

\title{Peano's Conception of a Single Infinite Cardinality}
\date{\today}

\begin{document}

\maketitle

\begin{abstract}

While Peano's negative attitude towards infinitesimals, in particular, \textit{geometric} infinitesimals, is widely documented, his conception of a \textit{single infinite cardinality} and, more generally, his views on the \textit{infinite}, are a lot less known. 
The paper reconstructs the evolution of Peano's ideas on these questions, and formulates several hypotheses about their underlying motivations. 
    
\end{abstract}

\begin{flushleft}

\section{Introduction}

Peano is widely known for his contribution to the creation of modern logical symbolism and to the axiomatisation of arithmetic, in particular, for the axioms which still bear his name.\footnote{For Peano's varied contributions to mathematical logic, see \citep{kennedy1980}, \citep{borga-freguglia-palladino1985}, \citep{grattan-guinness2000}, Ch. 5, and the 2021 special issue of \textit{Philosophia Scientiae} `Peano and His School' edited by P. Cant\`u and E. Luciano.} His proof of the impossibility of \textit{infinitesimal segments} in geometry, has also received the due attention;\footnote{\label{inf}Cf. \citep{peano1892a}. For a reconstruction of Peano's views about the infinitesimals, see \citep{ehrlich2006}, and the more recent \citep{freguglia2021}.} a lot less, if at all, explored are, on the contrary, Peano's views on the \textit{infinite},  an area which, precisely when Peano was most active mathematically, had experienced major and unprecedented breakthroughs, especially thanks to the work of, just to mention a few names, Cantor, Dedekind, Schr\"oder, Veronese and other mathematicians.\footnote{Cf., again, \citep{grattan-guinness2000} and \citep{ferreiros2010}.}

Recently, Peano's early conception of a \textit{single infinite cardinality} has been examined by \citep{mancosu2016}, which, besides showing how the conception locates itself in the history of the infinite, also casts `Peano's Principle' as a \textit{bona fide} cardinality principle putting pressure on the neo-logicist doctrine of the \textit{analyticity} of Hume's Principle.\footnote{\label{HP} Hume's Principle [HP], taken by Frege in \citep{frege1884} as the correct statement of `numerical identity' and, as a consequence, of the concept of number, is the universal closure of the following (second-order) principle: \[ \#x: F(x) = \#x: G(x) \leftrightarrow F \approx G \] which reads as follows: `the number of $F$s is equal to the number of $G$s if and only if $F$ is equinumerous with $G$ (ie., if $F$s and $G$s can be put in a one-to-one correspondence)'. The principle is, arguably, equivalent to Cantor's Principle, for which see later in the text, section \ref{galileo}. Mancosu has called the formalisation, in second-order logic, of Peano's conception of a single infinite cardinality `Peano's Principle': \[ \#x: F(x) = \#x: G(x) \leftrightarrow (\neg Fin(F) \wedge \neg Fin(G)) \vee ((Fin(F) \wedge Fin(G)) \wedge F \approx G) \]

\noindent
ie., `the number of $F$s is equal to the number of $G$s if and only if either $F$ and $G$ are infinite or, if finite, equinumerous', and discusses it in connection with the `good company problem' for HP, see \citep{mancosu2016}, Ch. 4.}

\citep{mancosu2016} also lists the different stages which Peano's conception of the infinite went through, and demonstrates that Peano's ideas became increasingly closer to, and from a certain point onwards, indistinguishable from, Cantor's, but it does not pin down the exact motivations behind the conception itself, nor does it explain why Peano, ultimately, became a full-blown supporter of Cantor's \textit{transfinite}.\footnote{\label{mancosu}Cf. \citep{mancosu2016}, pp. 165-6, in particular fn. 21.} 

The purpose of this paper is to fill in the gap, by providing, in a way as far as possible based on the extant textual sources, a more articulated account of the evolution of Peano's ideas, which might also shed light on the potential value and significance of Peano's early conception.   

The structure of the paper is as follows. In the next section, we make some further preliminary considerations. In \S 3, we examine the progressive development of Peano's ideas. In \S 4, we indicate two main potential motivations behind Peano's conception, and then proceed to make some final considerations in \S 5.

\section{Did Peano Have A Definite Conception of the Infinite?}\label{?}

To be sure, from reading Peano's terse statements on the subject, one wonders whether the Italian mathematician had \textit{any} determinate conception of the infinite at all; at times, it would rather seem that Peano's ideas are just insufficiently developed to be seen as advocating any specific view.

On some authors' judgement, Peano, while brilliant at spotting inaccuracies, and in improving on already established results, was poorly armed to devise (or just propose) a new theory, and rarely took enough pain in exposing in full detail, or in justifying, conceptions he happened to champion (or use). However, there is (some) evidence that this might not be the case with reference to the question of the infinite.\footnote{See, for instance, \citep{grattan-guinness2000}, p. 221, and Agazzi's introduction to \citep{borga-freguglia-palladino1985}, p. 7. For a different assessment of the significance of Peano's work, see \citep{rodriguez-consuegra1991}, especially the beginning of Ch. 3. }

To begin with, Peano actively engaged in the lively debates on the nascent theory of sets, and on infinitesimals, and relentlessly discussed aspects of these topics with Cantor, Russell, Frege, Veronese and other mathematicians.\footnote{Traces of such discussions may be found in, among other works, the many letters that Peano exchanged with his mathematical interlocutors, including Cantor (see later, section \ref{cantor}). The correspondence between Peano and other mathematicians may be found in \citep{peano2008}. } 
Moreover, as already mentioned, Peano forcefully opposed Veronese's introduction of the (geometric) infinitesimals, by producing a (purported) proof of their inconsistency,\footnote{See fn. \ref{inf} and section 4 of the present paper.} a fact which, on its own, already demonstrates that he was, at least, willing to go to great lengths to understand questions about the nature of the infinite, both philosophically and mathematically. 

Peano's ideas gradually became more and more akin to Cantor's, to the point that no traces of his earlier conception can be found in later work. What motivated such a noticeable change of mind? We believe that Peano became convinced of the intrinsic weakness of the alternatives to Cantor's conception, including his own, precisely as a consequence of his engagement in the foundational debates which were hosted by, among other journals, his \textit{Rivista di Matematica}. The correctness of this interpretation is further corroborated by the content of Peano's correspondence with Cantor in the year 1895, which, in our view, proved instrumental for Peano's `conversion' to the Cantorian approach. 

Therefore, especially in the decade 1890-1900, rather than showing laziness, or unwillingness to pursue the consequences of any specific conception, Peano looks very active in searching for the right approach to the issue of the infinite, by frequently interacting and discussing with his colleagues, and by testing and revising his position in a coherent manner. 

This said, even after inspecting much primary and secondary literature, as we have done in this paper, the full motivations behind Peano's early conception and behind the shift to the Cantorian conception, still appear, to some extent, mysterious, but the interpretative hypotheses we shall formulate in \S 4, will provide us with, at least, some clues.


\section{The Evolution of Peano's Conception}

In this section we follow closely the evolution of Peano's ideas between 1891 and the turn of the century. For this, we avail ourselves of three fundamental textual sources: 

\begin{enumerate}

\item Peano's own articles, in particular, \citep{peano1891}, where his early conception is first formulated, and \citep{peano1892a} (discussed in section \ref{CST});

\item Peano's correspondence with Cantor, consisting of six letters from Cantor to Peano dating to 1895;\footnote{Peano's letters to Cantor are no longer available, as confirmed, in a private email to the authors (23/02/2022), by Clara Silvia Roero, who has recently re-edited Peano's works and correspondence (cf., again, \citep{peano2008}).}
    
\item the first three volumes (and sections thereof) of his \textit{Formulaire de math\'ematiques}, published in rapid succession between 1895 and 1901 (\citep{peano1895}, then \citep{peano1897}, \citep{peano1898}, \citep{peano1899}, \citep{peano1901}), where both Peano's early and `Cantorian' conception may be found.

\end{enumerate}


\subsection{Peano 1891: `On the Concept of Number'}\label{1891}

In 1891, Peano published his famous article `Sul concetto di numero' (`On the Concept of Number') in the \textit{Rivista di matematica}, the journal he had founded the same year. 

Two years after \textit{Arithmetices principia, nova methodo exposita} (\citep{peano1889}), in which Peano first presented his axioms of the natural numbers, in this work, Peano, among other things, simplifies his notation, demonstrates that his famous five postulates of the natural numbers are mutually independent, and builds a more complete axiomatic system of numbers by also taking into account \textit{relative}, \textit{rational} and \textit{real} numbers.\footnote{Cf. \citep{borga-freguglia-palladino1985}, pp. 79-94, and \citep{rodriguez-consuegra1991}, Ch. 3.2.}

\citep{peano1891} is the first time Peano endeavours to deal explicitly with the infinite. In \S 9, he defines a function, `$num \ a$', whose domain consists of `classes' (denoted $a, b, c, ..., u,...$), and whose values are the `cardinalities' of these classes; in Peano's own words, $num \ a$ is `\textit{the number of elements of the class $a$}'\footnote{\label{cassina}[Con \textit{num a} intenderemo ``il numero degli individui della classe \textit{a}'' (p. 100).] Page numbers of \citep{peano1891} are those of the third volume of the Cassina edition of Peano's works: \textit{Opere scelte}, vol. 3, Cremonese, Roma, 1959. The English translations of Peano's quotes are all ours.}. 

Now, if $a$ is a finite class, then $num \ a$ is just the (finite) number of its elements, ie., a \textit{natural number} $n$. But then, Peano states that $num \ a$ is not always a natural number, since the set of natural numbers does not include `zero' and `infinity'.\footnote{[Data una classe $a$ non sempre $num(a)$ è un $N$, poich\'e $N$ non comprende n\'e lo zero, n\'e l'infinito (p. 101).]} This is the first time Peano mentions the possibility that a class $a$ be \textit{empty}, or \textit{infinite} (that is, that $num \ a=\infty$).\footnote{[Il segno $num$ è un segno d’operazione che ad ogni classe fa corrispondere o un $N$, o lo $0$, o l'$\infty$ (p. 102).]} 

We wish to say something more substantial about `$\infty$'. Peano seems to take it to be a \textit{bona fide} `infinite quantity', which can be manipulated like any other (finite) quantity, as is clear from the propositions 3 and 4 in \S 9:

\begin{quote}

3. If $a$ and $b$ are two non-empty and finite classes having no element in common, then the number of elements of the set of the two classes $a$ and $b$ is equal to the sum of the number of $a$s and $b$s.\footnote{[Essendo $a$ e $b$ due classi non nulle e finite, non aventi alcun individuo comune, allora il numero degli individui appartenenti all’insieme delle due classi $a$ e $b$ vale la somma dei numeri degli $a$ e dei $b$ (p. 102).]}
    
\end{quote}

\noindent
In 3. above, Peano is stating, in modern set-theoretic notation, that if two sets $a$ and $b$ are disjoint, then the number of the elements of $a$ \textit{and} $b$ is equal to the number of the elements of $a \cup b$. Peano, then, notes that the proposition holds even if one of the two classes, and even both, contain \textit{infinite} elements; but now, we have, as a consequence, that:

\begin{equation}
    \label{prop3_1}
    x + \infty = \infty + x = \infty
\end{equation}

and 

\begin{equation}
   \label{prop3_2}
   \infty + \infty = \infty. 
\end{equation}

\noindent
Immediately afterwards, he says:

\begin{quote}

4. If the classes $a$ and $b$ are such that the second is contained in the first, and the class $b$ is non-empty, and is not equal to $a$, and if the number of $a$s is finite, then the number of $b$s is also finite, and is less than the number of $a$s.\footnote{[Se delle classi $a$ e $b$ la seconda è contenuta nella prima, e la classe $b$ non è nulla, e non è eguale ad $a$, e se il numero degli a è finito, allora anche il numero dei $b$ è finito, ed è minore del numero degli $a$ (p. 101).]}

\end{quote}

\noindent
Then he notes: `this proposition ceases to be valid if $num \ a=\infty$'.\footnote{[Questa proposizione cessa di esistere se $num \ a=\infty$ (p. 101).]} Propositions 3 and 4 also demonstrate that `$\infty$' is taken by Peano to be different from Cantor's `$\omega$', since, as is known, by Cantor's conception, $\omega+n \neq n+\omega$.  Indeed, on the grounds of the arithmetical laws outlined in \eqref{prop3_1} and \eqref{prop3_2}, one could be tempted to view `$\infty$' as being equivalent to `$\aleph_0$'. But this would be a hasty conclusion. In section \ref{CST}, we shall see that Peano's own intepretation of his `infinitary numbers' does not automatically sanction the equivalence between his `$\infty$' and `$\aleph_0$'. 

One further comment is in order. Proposition 4 states that infinite classes could contain infinite \textit{proper} subclasses. This implies, among other things, that `standard' part-whole relationships cease to be valid in the infinite, a fact which, in turn, may have had some bearings on the development of Peano's ideas, as we shall see in the next section.

To sum up, \citep{peano1891}'s conception clearly is non-Cantorian, as it fosters the existence of just one infinite `cardinal number'. So the question arises of whether and to what extent Peano had assimilated Cantor's theory at the time. On this issue, it should be preliminarily observed, on the one hand, that only one year later, in \citep{peano1892a}, Peano exhibits a more than decent understanding of Cantor's transfinite. On the other, as we shall see in a moment, his 1895 correspondence with Cantor, in particular, the questions he asked his German colleague, reveal that, at the time of the composition of \citep{peano1891}, Peano still felt the need to understand Cantor's theory in more depth. 

\subsection{The Cantor-Peano Correspondence (1895)}\label{cantor}

Peano's early conception of the infinite, which we have just commented upon, is gradually modified, and then definitively replaced, in 1899, by Cantor's theory of the transfinite. What, ultimately, led Peano to change his mind remains an open question. Although the evidence is insufficient, we conjecture that the correspondence between Peano and Cantor, in the year 1895, proved instrumental in that respect.  

In what follows, we review bits of the discussion Peano entertained with Cantor that are relevant to our purposes; as already noted, one side of the correspondence (from Peano to Cantor) is not extant, so Peano's comments and answers can only be, very approximately, deduced from the responses of his German colleague.\footnote{Cantor's letters to Peano examined in the present work are all collected in the Meschkowski edition of Cantor's letters, \citep{cantor1990}, pp. 359ff. (n. 143-147). An overview of the contents of the Cantor-Peano correspondence is in \citep{kennedy1980}, pp. 87-90.}  


The letters cover the following subjects:

\begin{enumerate}
    \item Cantor's articles that Peano intends to publish in his \textit{Rivista Matematica};
    
    \item Cantor's distaste with Veronese's work on the infinite, which posited the existence of infinitesimals;
    
    \item Peano's request for explanations about Cantor's elusive definition of `cardinal number' in Cantor's works.
\end{enumerate}

In particular, (2) and (3) are fundamental, in our view, for the evolution of Peano's ideas on the infinite. Let's examine them in more detail.

As far as (2) is concerned, in the two letters of, respectively, July 27 and July 28, 1895, the discussion focuses on the conflict between Cantor's and Veronese's theories. On the German mathematician's view, what Veronese called `ordered groups' were just a plagiarism of Cantor's `simply ordered sets'. But then, Cantor points out what he thinks is the mistake that Veronese has made in trying to automatically extend the arithmetic of natural numbers to infinite numbers. As Cantor observes, in two passages of the July 27 letter: 

\begin{quote}

But if it is correct that his $\infty_1 = \omega+^{*}\omega$, his assertion that: \[2 \cdot \infty_1 = \infty_1 \cdot 2 \]

must be \emph{incorrect}! [...] Anyway, his [Veronese's, \textit{our note}] `infinite numbers' seem tenable to me only if they are identified with some of my `transfinite order-types'. In this case, however, they lack the law of commutativity for addition and multiplication (in general) on which he [Veronese, \textit{our note}] lays so much stress.\footnote{[Jedenfalls scheinen mir seine ”unendlichen Zahlen” \textit{nur dann haltbar, wenn sie mit gewißen von meinen "transfiniten
Ordnungstypen” identificirt werden}. In diesem Falle fehlt ihnen aber das
Gesetz der Commutabilität bei der Addition und Multiplication (im Allgemeinen), \textit{worauf er solchen Nachdruck legt} (\citep{cantor1990}, p. 360)]. The English translation is ours.}
    
\end{quote}

\noindent
We do not know the content of Peano's answer; in any case, on Cantor's own impulse (28 July 1895 letter), the 27 July letter was published by Peano in the \textit{Rivista di Matematica}.\footnote{In the August issue of the journal.} So, it is imaginable that Cantor had managed to convince Peano of the correctness of his arguments.

As far as (3) is concerned, in the subsequent letters, the conversation between the two mathematicians turns to consider Peano's qualms about Cantor's definition of `cardinal number'. This can be safely deduced from Cantor's September 14, 1895 letter, in which Cantor tries to clarify some notions contained in the section 5 of \citep{cantor1895}. From what Cantor says, it seems quite clear that Peano was unsure about Cantor's notion of `finite cardinal number' and about how Cantor introduced the induction principle, and asked his colleague for clarifications.

In the 21 September 1895 letter, Cantor quotes a sentence written by Peano in his response to Cantor's previous letter: `\textit{Where can one find the definition of finite cardinal numbers?}',\footnote{[Où est ce que l'on trouve la definition der endlichen Cardinalzahlen? (\citep{cantor1990}, p. 365)].} a clear sign that Cantor's first letter was not entirely clarificatory. 

Therefore, what is likely is that, as a consequence of the clarifications Cantor gave Peano, and of Cantor's rebuttal of Veronese's theory, Peano started considering Cantor's conception of the transfinite as being the only correct one, and Cantor's theory of cardinal numbers as pretty much definitive. If in 1891 he had chosen not to fully adhere to Cantor's theory - or simply not to explore it in depth - now he had further material at hand, including Cantor's clarificatory statements, which could re-orient his views, a fact which would soon be reflected, as we will see, in the \textit{Formulaire}, where the theory of the transfinite gradually takes primacy.

\subsection{The Infinite in the \textit{Formulaire} (1895 to 1901)}\label{form}

We finally turn to survey, very briefly, the modifications of Peano's conception that can be found in the volumes (editions) of Peano's \textit{Formulaire de math\'ematiques}.\footnote{The volumes of the \textit{Formulaire} all appeared as supplements to issues of the \textit{Rivista di Matematica} between 1895 and 1908. For an overview of the plan and evolution of Peano's \textit{Formulaire}, or \textit{Formulario} (in latino sine flexione), see the very exhaustive \citep{cassina1955}.} As said at the beginning, these have already been briefly taken into account by Paolo Mancosu in \citep{mancosu2016};\footnote{See fn. \ref{mancosu}.} what follows aims to expand on, and complement, Mancosu's account. 

The evolution of Peano's ideas on the infinite in the \textit{Formulaire} spans a period of 4 years, and ends in 1899, when all traces of Peano's early conception disappear, and Cantor's transfinite ultimately takes over.   

The first volume (\citep{peano1895}) saw the collaboration of other Italian mathematicians, such as G. Vailati, F. Castellano, G. Vivanti, R. Bettazzi. Of particular importance for our purposes are sections V and VI, written, respectively, by Peano and Giulio Vivanti; the latter, among other things, had taken part in the debate with Bettazzi about the infinitesimals raging in the \textit{Rivista} in the years 1891-1892.\footnote{The debate is accurately reconstructed in \citep{ehrlich2006}, pp. 75-101.} 

In section V of \citep{peano1895}, entitled `Classes de nombres', we find again the mentioning of just one infinite cardinality, `$\infty$', which is now explicitly defined as one of the possible values of $num \ u$, where $u$ is a class, as follows:\footnote{In the text of formulas (4)-(6) below, for the sake of simplicity, we have chosen not to adopt Peano's original notation, but the modern one.}

\[ 4. \ num \ u=\infty \leftrightarrow num \ u \notin \mathbb{N}_{0} \]

\noindent
Then Peano explains that $num \ u$ can take only three values:

\[ 5. \ num \ u= n \in \mathbb{N} \vee 0 \vee \infty\]

\noindent
Moreover, `$\infty$', in definition 6, is, again, characterised as enjoying \textit{commutativity}:

\[6. \ a \in \mathbb{N}_0, \quad a+\infty=\infty + a = \infty+\infty=\infty \quad and \quad  a<\infty\]

\noindent
On the other hand, section VI, due to Vivanti, not to Peano, introduces, and uses, the symbols `\textit{Nc}' and `\textit{Ntransf}', denoting, respectively, Cantor's cardinals and ordinals, and deals with the corresponding set-theoretic notions. Therefore, in 1895, Peano did not entertain ideas much different from those appeared in \citep{peano1891}, ie. he did not subscribe to Cantor's views, while \citep{peano1895}'s section VI, composed by Vivanti, which mentions Cantor's cardinals and ordinals, has no immediate connection with the one composed by Peano. This is consistent with our hypothesis that the correspondence with Cantor may have led Peano to change his mind, as the text of \citep{peano1895} must have been produced well before 1895.

The second volume of the \textit{Formulaire} consists of three parts, each published one year apart from the other, following the editorial plan below:

\begin{itemize}
\item \S 1, \textit{Logique mathématique} \citep{peano1897};
\item \S 2, \textit{Arithmétique} \citep{peano1898}
\item \S 3 (untitled) \citep{peano1899}

\end{itemize}

\citep{peano1898} contains the first manifestation of Cantor's theory in Peano's writings.  Peano says: 

\begin{quote}

Here we define another idea, indicated by the sign $Nc'$ similar to the previous one, indicated by $num$, but not identical. This definition is expressed by the only signs of logic (\S1); it is therefore independent of the primitive ideas $\mathbb{N}_0, +, 0$. We could start arithmetic here. P211 expresses the coincidence of the signs $num$ and $Nc'$ when one deals with finite classes. But we will have for example: $num \ \mathbb{N} = num \ \mathbb{R^+}$, because the classes $\mathbb{N}$ and $\mathbb{R^+}$ are both infinite;\footnote{In the original text, Peano uses the symbols N and Q, to denote, respectively, the set of the \textit{natural} and of the positive \textit{real} numbers; in our translation, for the sake of simplicity, we have replaced N and Q with the modern $\mathbb{N}$ and $\mathbb{R^+}$.} but $Nc' \ \mathbb{N}< Nc' \ \mathbb{R^+}$, because the power of $ \mathbb{N}$ is less than the power of $ \mathbb{R^+}$. Mr. Cantor indicated the power of $a$ by $\overline{\overline{a}}$ (See RdM. A. 1895 p.130), a notation that cannot be adopted in the Formulaire. M. Vivanti in $F_1$ VII §2 P1 [section VI of \citep{peano1895}, \textit{our note}] has replaced it with $Nc' \ a$ ``the cardinal number of $a$.''.\footnote{[Nous définissons ici une autre idée, indiquée par le signe ``$Nc'$'' semblable à la précédente, indiqueée par ``num'', mais non identique. Cette définition est exprimée par les seules signes de logique (\S 1); elle est donc indépendante des idées primitives $\mathbb{N}_0, +, 0$. On pourrait commencer ici l'Arithmétique. La P211 exprime la coincidence des signes $num$ et $Nc'$ lorsqu'il s'agit de classes finies. Mais on aura par exemple: $num \ \mathbb{N} = num \ \mathbb{R^+}$, car les classes $\mathbb{N}$ et $\mathbb{R^+}$ sont toutes les deux infinies; mais $Nc'\ \mathbb{N}<Nc' \ \mathbb{R^+}$, car la puissance des $\mathbb{N}$ est plus petite que la puissance des $\mathbb{R^+}$.  M. Cantor a indiqué la puissance de $a$ par $\overline{\overline{a}}$ (Cfr. RdM. a. 1895 p.130), notation qu'on ne peut pas adopter dans le Formulaire. M. Vivanti dans $F_1$ VII §2 P1 l'a remplacée par $Nc' \ a$ "le nombre cardinal des $a$''. (\citep{peano1898}, p. 39)].}
    
\end{quote}

\noindent
From the quote above, it seems clear that Peano is beginning to change his mind, and gradually converting to Cantor's theory, although he still keeps mentioning both conceptions: $Nc'$, denoting Cantor's cardinal numbers, now appears alongside $num$ as cardinality assignment to an infinite class.

\citep{peano1899} marks one further transformation of Peano's conception. In that work, for the first time, Peano `merges' the concepts of $num$ and $Nc'$ into a unique, and new, concept denoted by `$Num$'. $Num$'s are now nothing but Cantor's transfinite cardinalities, an `expansion' of Peano's former $num$ subtly revealed by the change of notation. From \citep{peano1899} onwards, thus, Peano seems to fully adhere to Cantor's theory. 

In the \textit{Formulaire}'s third volume, \citep{peano1901}, in the $Num$ section, Peano, finally, explains:

\begin{quote}

$Num'Cls$ means ``the number of a class''. These numbers coincide with the [elements of] $\mathbb{N}_0$ for the finite classes; G. Cantor calls them ``cardinal numbers''. In F 1895 the symbol “$Nc$” was introduced to represent them [again, a reference to Vivanti's section VI of \citep{peano1895}].\footnote{[$Num'Cls$ signifie ``le nombre d'une classe''. Ces nombres coïncident avec les $N_0$ pour les classes finies; G. Cantor les appelle ``nombres cardinaux''. Dans F1895 on a introduit le symbole ``$Nc$'' pour les représenter. (\citep{peano1901}, p. 70)]}
    
\end{quote}

Thus, already by 1901, and then in the subsequent volumes of the \textit{Formulaire} (1903, 1908), any significant distinction between Peano's conception and Cantor's disappears, and Peano's `single-cardinality' conception of the infinite is a distant memory.


\section{Interpretations of Peano's Conception}

In this section, we turn to examine the issue of what may have motivated Peano's conception of a single infinite cardinality and present further textual sources which might help us understand better the evolution of Peano's ideas. 

We will follow two main interpretative routes: (1) the first will highlight Peano's hesitancy to either adopt the Cantorian conception of infinite cardinalities as opposed to the one based on the Part-Whole Principle; our underlying presupposition is that, like Galileo before him, Peano was not fully convinced of the correctness of either method; (2) the second will point out an `incomplete' understanding, or misuse, of Cantor's theory of ordinals and cardinals, which may have misled Peano to believe that his theory could be `reduced' to Cantor's or, in alternative, that Cantor's theory of the transfinite could be bent to suit his own purposes. We examine, respectively, (1) in \S 4.1, and (2) in \S 4.2.

\subsection{A Galilean View?}\label{galileo}

In his \textit{Dialogues Concerning Two New Sciences} (1638), Galileo came to express skepticism about the possibility of carrying out size comparisons with regard to infinite `collections' of objects. 

In the work's \textit{First Day}, Galileo's spokesperson, Salviati, proposes to compare the magnitude of a collection of objects with that of another collection by checking that each object in the former collection is \textit{uniquely} matched by one in the latter. His case study consisted in comparing square with natural numbers (in modern parlance, we would say the set of square numbers ($\mathbb{S}$), and that of natural numbers ($\mathbb{N}$)), and the mapping function that he is using is the one that associates each natural number with its square:  $n \mapsto n^{2}$:

\[ 0 \leftrightarrow 0\]
\[ 1 \leftrightarrow 1\]
\[ 2 \leftrightarrow 4\]
\[ 3 \leftrightarrow 9\]
\[ .. \leftrightarrow ..\]

\noindent
Salviati and his interlocutor, the learned listener Sagredo, agree that, using this method, it must be concluded that there are \textit{as many} numbers \textit{as} squares (that is, that $\mathfrak{s}(\mathbb{N})=\mathfrak{s}(\mathbb{S})$, where $\mathfrak{s}(X)$ denotes the `size of $X$'). But now the dire puzzle unfolds in front of their eyes: it seems to be an `established fact', in any case, commonsense knowledge, that there are \textit{more} natural \textit{than} square numbers, since $\mathbb{S}$ has lots of gaps `in between', that $\mathbb{N}$ does not have; this is natural, since $\mathbb{S}$ is a \textit{proper part} of $\mathbb{N}$ (ie., $\mathbb{S}\subset \mathbb{N}$). So, relying on intuitions referring to the `density' of $\mathbb{S}$ in $\mathbb{N}$, Salviati suggests that one should, in fact, more correctly, conclude that $\mathfrak{s}(\mathbb{N})>\mathfrak{s}(\mathbb{S})$.  This is what has come to be known as Galileo's Paradox. 

Now, asks Salviati, which of the two horns of the dilemma is correct? Neither, apparently, as Galileo ultimately points out, since:

\begin{quote}

This is one of the difficulties which arise when we attempt, with our finite minds, to discuss the infinite, assigning to it those properties which we give to the finite and limited; but this I think is wrong, for we cannot speak of infinite quantities as being the one greater or less than or equal to another. (\citep{galileo1638}, p. 31)
    
\end{quote}

\noindent
Now, Galileo does not go as far as to say that, but, since one cannot speak of magnitudes `greater than', `less then', or `equal to' one another in the \textit{infinite}, it would, in principle, be fully legitimate and consequential to think that, to all purposes and intents, there is but \textit{one} infinite cardinality. Thus, what we would have at hand, in \citep{galileo1638}, would be the enunciation of a `single-cardinality' conception of the infinite dictated by the impossibility to come to terms with a conflict between two ways of counting in the infinite, one based on `bijections' and one on `density'. The former method is formally enshrined in what would, much later, become the central pillar of set theory, that is: 

\vspace{11pt}

\noindent
\textbf{Cantor's Principle (CP)}. \textit{Given two sets $A$ and $B$, if there exists $f: A \rightarrow B$ which is 1-1 and onto (ie., a bijection between $A$ and $B$), then $\mathfrak{s}(A)=\mathfrak{s}(B)$.}\footnote{CP may be seen as a `variant' of Hume's Principle mentioned in \S 1 (see fn. \ref{HP}).} 

\vspace{11pt}

\noindent
The latter obeys what was already known in antiquity as Euclid's Axiom (Common Notion V of the \textit{Elements}), that is, the: 

\vspace{11pt}

\noindent
\textbf{Part-Whole Principle (PWP)}. \textit{Given two sets $A$ and $B$, if $A \subset B$, then $\mathfrak{s}(A)<\mathfrak{s}(B)$.}

\vspace{11pt}

\noindent
In our era, the `set-theoretic era', as it were, one tends to think quite immediately that CP is the correct way to assign cardinalities to infinite collections, and that, as a consequence, Galileo's Paradox is no paradox at all, but rather helps exhibit the hallmark of the infinite itself, as captured by the notion of:

\vspace{11pt}

\noindent
\textbf{Dedekind-Infiniteness}. A set is said to be \textit{infinite} if and only if it can be put in a \textit{one-to-one correspondence} with a proper part of itself; otherwise, it is \textit{finite}.\footnote{The notion was first formulated in \citep{dedekind1888}, V, 64.} 

\vspace{11pt}

\noindent
But the correctness, and indispensability, of PWP, as is known, was the prevailing view for many centuries, and still lay almost uncontested in Galileo's times.\footnote{For a careful excursus of the history of several cardinality principles, we refer the reader to \citep{mancosu2009}.} Moreover, oscillations between compliance with CP or PWP could even be found, centuries later, in avowed supporters of the `bijection method' like Bolzano.\footnote{Cf. \citep{mancosu2016}, pp. 130-1.}

Now, one could speculate that: 1) Peano shared Galileo's concerns about the possibility of articulating a theory of \textit{distinct} infinite cardinalities, and that: 2) he retreated to a more sober `one-cardinality' conception precisely because he was hesitant to choose among one of the two (known) methods of counting in the infinite.

In particular, Peano's 1891 conception might have reflected dissatisfaction both with a purely Cantorian \textit{and} with the Euclidean point of view.  

This is (indirectly) proved by a remark made by Peano in a footnote of his \citep{peano1891}, concerning Rodolfo Bettazzi's use of the `bijection method' to characterise the concept of `cardinal number' in \citep{bettazzi1887}.\footnote{The passage is cited (and translated to English) by \citep{mancosu2016}, p. 164-5.} Apparently, Bettazzi thought that, if two sets $A$ and $B$ could be put in a one-to-one correspondence, then any correspondence between them would inevitably be one-to-one. For the sake of our argument, let us state this as:

\vspace{11pt}

\noindent
\textbf{Bettazzi's Principle (BP)}. \textit{Given two sets $A$ and $B$, if there exists a bijective function between $A$ and $B$, then there is no function between $A$ and $B$ which is injective, but not surjective.}

\vspace{11pt}

\noindent
Peano correctly objects to Bettazzi that one can biject infinite sets with their own infinite parts, that is, one can have \textit{injections} of sets with themselves which, evidently, are not \textit{surjective}. 



Considerations such as these could have led Peano to doubt that both CP and PWP could be a satisfactory characterisation of the notion of `infinite cardinality', and may ultimately have convinced him that all infinite classes had the same, unique \textit{cardinality}. This can further be deduced from two facts. One the one hand, if Peano had accepted PWP as the correct conception of infinite cardinalities, then he would have had to admit of different orders of infinity, that, for instance, of a subset of $\mathbb{N}$ and of $\mathbb{N}$ itself. But this he doesn't do either in \citep{peano1891} or nowhere else, for that matter. On the other, as we have already seen in section 3.3, Peano does not mention the existence of distinct infinite cardinalities as based on the `bijection method', and, as late as \citep{peano1895}, he still maintains that $num(\mathbb{N})=num(\mathbb{R^+})$. 

There is still one further way in which Peano's position could be viewed as Galilean. If one just takes into account $\mathbb{N}$ and its subsets, one may be perfectly entitled to conclude that there is but \textit{one} infinite cardinality (maybe recognition of this fact already underlay Galileo's own thinking). Now, with the exception of later examples he makes of infinite cardinalities where he also takes into account, e.g., $\mathbb{R}$, Peano's considerations (especially those in \citep{peano1891}) mostly refer to sets and subsets of natural numbers. 

In sum, assuming that our `Galilean' interpretation of Peano's conception is correct, one might conclude that either Peano felt that there was a tension between CP and PWP in the infinite which, ultimately, implied that there were no distinct  infinite cardinalities, or he thought that CP, restricted to $\mathbb{N}$ (and subsets thereof), that is, to the realm of \textit{arithmetic} in which he was especially interested, licensed the view that there is just \textit{one} infinite cardinality.

\subsection{(Geometric) Infinitesimals and Cantorian Set Theory}\label{CST}

We have expressed many times the view that Peano may have progressively embraced Cantor's transfinite because he finally came to grips with it after a resolutive exchange of letters with Cantor himself which clarified in full the details of Cantor's theory. There are, however, passages of Peano's earlier works which do not lend support to this thesis. In particular, we are referring here to the article on infinitesimals, \citep{peano1892a}, where Peano aims to improve on Cantor's earlier proof that infinitesimal segments are \textit{inconsistent}.\footnote{A detailed analysis of Cantor's proof (for which see \citep{cantor1887-88}) is carried out in \citep{ehrlich2006}, pp. 27-51.}  In the article, Peano shows a good command of Cantorian set theory, but seems to want to use it, very originally, in a way compatible with his `single-cardinality' conception of the infinite. Let us see how in more detail. 

The key notions in that work are those of `bounded' and `infinitesimal segment'. Given the half-line with origin in $o$, a `bounded' segment is a segment $op$ with origin in $o$ and end in $p$. Peano also denotes a bounded segment with $u$, or with $\overline{ou}$. An `infinitesimal segment' is a $\overline{ou}$ lying inside a bounded segment $v$, and denoted $u=v/\infty$. The latter's full definition reads as follows:

\begin{quote}
    We say that the \textit{segment} $u$ is infinitesimal with respect to $v$ and we write $u \ \epsilon \ v/\infty$, if every multiple of $u$ is less than $v$ [...].\footnote{[Dicesi che il segmento $u$ è infinitesimo rispetto al segmento
$v$, e scriveremo $u \ \epsilon \ v/\infty$, se ogni multiplo di $u$ \`e minore di $v$ [...] (p. 113)] Note the use of Peano's $\epsilon$ symbol, which means: `is'. Page numbers for \citep{peano1892a} are, again, those of Cassina's edition, for which see fn. \ref{cassina}. } 
\end{quote}

\noindent
Then Peano proceeds to define the multiple of infinite order of $u$, `$\infty u$'. He says:

\begin{quote}
    We shall posit:
    
    \[\infty u = \bigcup \mathbb{N}u \]
    
that is, we call multiple of $u$ of infinite order the set of points which either lie on some segment $u$, $2u$, $3u$, $...$, or the upper bound of the multiples of $u$.\footnote{[Porremo \[\infty u = \bigcup \mathbb{N}u \] cioè chiamiamo multiplo d'ordine infinito di $u$ l'insieme dei punti che stanno sopra qualcuno dei segmenti $u$, $2u$, $3u$, $...$, o il limite superiore dei multipli di $u$. (p. 113)] For the sake of simplicity, in the formula, we have chosen to replace Peano's original notation: $\cup'Nu$ with $\bigcup \mathbb{N}u$.}

\end{quote}

\noindent
Now, Peano asserts that, by the definition of infinitesimal, if $u$ is infinitesimal, also `$\infty u$' must lie \textit{inside} $v$, and so must be all other infinitary multiples of $u$. In a crucial passage of \citep{peano1892a}, he says:

\begin{quote}

We can add $\infty$ to itself, thus obtaining 
$2\infty u$, and, generally, we can form all multiples of $\infty u$; we can multiply $\infty u$ by $\infty$, and obtain $\infty^{2} u$ and so on. But all these various segments, which one obtains multiplying  $u$ by Cantor's transfinite numbers are all equal to one another [...].\footnote{[Possiamo sommare $\infty u$ con s\'e
stesso, ottenendo così $2\infty u$, ed in generale possiamo formare tutti i multipli di $\infty u$; possiamo moltiplicare $\infty u$ per $\infty$, ed ottenere
$\infty^{2} u$ e così via. Ma tutti questi varii segmenti, che si ottengono moltiplicando
$u$ pei numeri transfiniti di Cantor sono eguali fra loro, [...]] (p. 113)}
\end{quote}

\noindent
As a consequence, infinitesimal segments, if such things existed, would blatantly violate the properties of bounded segments, which require that, for instance, a segment of length $2 \infty u$, twice as much as $\infty u$, be \textit{greater than} $\infty u$. Hence, concludes Peano, infinitesimal segments are inconsistent with `standard' representations of the geometric space. 

Leaving aside entirely the issue of the correctness of Peano's argument, what is most relevant for our purposes is to try and understand Peano's reference, in it, to `Cantor's transfinite numbers'. The mentioning of such numbers is indeed puzzling, as, by Cantor's own theory, as we know, $\omega+1>\omega$, hence, $\omega+1 \cdot u$ should be greater than $\omega \cdot u$. 

Indeed, Giuseppe Veronese, the main target of Peano's (and Cantor's) polemic against geometric infinitesimals, was able to react very quickly to Peano's argument. In \citep{veronese1892}, which is meant to be a response to \citep{peano1892a}, the Italian mathematician correctly diagnoses what seems to be the trouble with Peano's proof: 

\begin{quote}
    But these equalities [$\infty u = 2\infty u$ = $\infty^{2}u = ...$, \textit{our note}] do not depend on the properties of Mr. Cantor's transfinite numbers, for which hold: $\omega+1>\omega$, $2\omega>\omega$, etc., but precisely on considering $\infty u$ as \textit{unlimited} (our italics).\footnote{\citep{veronese1892}, p. 74.}
\end{quote}

\noindent
Veronese is, here, drawing attention to the fact that it is the very definitions of $\infty u$, $2\infty u$, etc., not \textit{transfinite arithmetic}, which allow one to deduce that these quantities are all equal. 
Similar remarks were made by Giulio Vivanti, who, in his \citep{vivanti1895}, explained why the argument was doomed to failure. Vivanti argued that Peano viewed $(\infty+1)u$ as the upper bound of $(n+1)u$, but the points of $(\infty+1)u$ must lie either on one of the finite multiples of $u$, or be $(\infty+1)u$ itself. So, $(\infty+1)u$ is precisely the same the segment as $\infty u$; hence, $(\infty+1)u$ and $\infty u$ must be \textit{equal}. But, if it is true that Cantor's $\omega$ is the limit of $n$, $\omega+1$ is no \textit{limit number (ordinal)}, and, as already pointed out by Veronese, $\omega+1>\omega$, and, consequently, $(\omega+1)u>\omega u$.\footnote{\citep{vivanti1895}, p. 69.} 

So, Peano's reliance upon `Cantor's transfinite numbers' in his proof does not serve well his purpose of showing that geometric infinitesimals are inconsistent. In order to rescue the force of the argument, \citep{freguglia2021} conjectures that the reason why Peano holds the correctness of the equation:

\[ \infty u = 2 \infty u = \infty^{2} u = ... \]

\noindent
is the fact that Peano is, intentionally or not, `assimilating' `$\infty$' to `$\aleph_0$'.\footnote{Cf. \citep{freguglia2021}, p. 152ff. The author says: `Peano explicitly assimilates $\infty$ to $\aleph_0$ [...]' and then explains ibid., fn. 14: `[i]n the sense that it has the same arithmetic behaviour.'} This is because, Freguglia argues, while, as said, $\omega \neq \omega+1$, on the contrary, $\aleph_0=\aleph_0+1$ and, in general, if $\kappa$ is a transfinite \textit{cardinal} number, and $n$ a natural number:

\[ \kappa+n = \kappa \cdot n = \kappa^{n} = \kappa \]

\noindent
Freguglia's interpretation has some merit, but seems to be missing, at large, the mark. In particular, it seems to be rather unjustified to force this interpretation on Peano's reference to `Cantor's transfinite numbers'; the numbers Peano is referring to here are, it seems to us, and as understood by Veronese and Vivanti, Cantor's \textit{transfinite ordinals} ($\omega, \omega+1, ..., \omega+\omega$, ...).

However, following the spirit, not the letter, of Freguglia's intepretation, one could say that Peano may have taken his numbers ($\infty u$, $\infty^{2}u$, ...) to be equivalent to $\aleph_0$ in the sense that, although, like Cantor's ordinals they were linearly ordered by `$<$', and, thus, differed from each other in `magnitude', they had, nonetheless, the same \textit{cardinality}, in particular, they were all \textit{countable}. Now, the reasoning goes, since what counts as \textit{length} of a geometric segment is a measure expressed by a cardinal number, one must conclude that segments of length $\infty u$, $\infty^{2} u$ are all \textit{equal}. But there's a problem with this interpretation: Cantor's ordinals may be \textit{uncountable}, that is, they may have different cardinalities, and Peano might have been aware of this. 

We wish to propose an alternative interpretation of Peano's argument, which aligns more with the content of Vivanti's and Veronese's comments above. This avers that Peano viewed his `infinite numbers' as being all equal to each other not because he thought that `$\infty$' was equivalent to `$\aleph_0$', but because he thought that there was just \textit{one}, maximal, unsurpassable, infinite cardinality, `$\infty$'. By this interpretation, when Peano is referring to `Cantor's numbers' in the aforementioned passage, he is really likening his own numbers to Cantor's ordinals (as supposed by Veronese and Vivanti), but he is also committing himself to the view that those numbers must all have the same cardinality, because, by his own conception, there exists \textit{only} one, `$\infty$', which cannot be transcended. So, in a sense, Peano's reference to `Cantor's numbers' can be seen as an attempt to make Cantor's conception suit (to some extent, even merge it) Peano's own, `single-cardinality', conception. 

Our interpretation entails that Peano had developed a carefully thought of, and original, conception of the infinite in the years 1891-1895, which he believed could be employed in practical mathematical contexts (such as that discussed in \citep{peano1892a}). By this very interpretation, the fact that he, subsequently, abandoned this conception does not mean that he thought it to be faulty, `half scarce made up' in any respect, but only that he thought, as we have have already said many times, Cantor's theory to be more general, far-reaching and, ultimately, fruitful than his own.  

\section{Concluding Remarks and Ways Forward}

We have reviewed the development of Peano's conception of the infinite up to the 1901 publication of the third volume of his \textit{Formulario}. We have seen that, already by the publication of the third section of the second volume in 1899, Peano's conception had aligned with the, then already `meainstream', set-theoretic conception of infinite cardinalities. 

In section \ref{cantor}, we have conjectured that Cantor may have played a role in persuading Peano to abandon his earlier view. As far as the latter is concerned, we have shown that, in his \citep{peano1892a}, Peano showed to have perfectly clear in his mind the consequences of his own conception, the way it could be properly used in `practical' mathematical contexts, and how it related to Cantor's transfinite. Moreover, we have seen that, up to a point, Peano might have been doubtful about the efficacy of the `bijection method', as this conflicted with PWP, and we have conjectured that these doubts might have led him to formulate the conception of a single infinite cardinality. 

By way of conclusion, we wish to draw attention to some recent developments, which could be, if loosely, connected to Peano's views. 

The debate on alternatives to Cantor's conception of cardinality has, lately, flourished again, in ways which potentially suggest a re-thinking of the primacy of CP. The concrete way this is usually done is through re-legitimising the use of PWP.\footnote{For recent, alternative conceptions of infinite cardinality based on PWP see Sergeyev's theory of \textit{grossone} (cf. \citep{sergeyev2017}), and Benci et al.'s theory of \textit{numerosities} (\citep{benci-dinasso-forti2006}).}

For another development, one could attempt to narrow down the extent of Cantor's transfinite by introducing suitable restrictions on generating principles for sets, as imagined, for instance, by Vop\v{e}nka's `alternative set theory' with just two infinite cardinalities, the \textit{countable} and the \textit{uncountable}, Randall Holmes' `pocket set theory', or Meir Buzaglo's cardinality principles discussed in \citep{buzaglo2002}.\footnote{For Vop\v{e}nka's `alternative set theory', see \citep{vopenka1979}. For `Buzaglo's Principles', see \citep{buzaglo2002}, pp. 43-4, 49-50. `Pocket set theory' is outlined in \citep{rholmes2017}, pp. 62ff.}

Just to make a concrete example, let us review one of Buzaglo's `principles'. This is based on defining the `$\nless$' (`not smaller than') relationship between the sizes of sets, in turn the negation of a special `$<$' relationship defined as follows:\footnote{The `$|...|$' sign, in the definitions below, should be taken to denote a generic concept of `size of a set', not the set-theoretic one.} 

\begin{definition}[Smaller Than]

$|A|$ is smaller than $|B|$ if and only if: (i) there is a function $f$ from $A$ to $B$; (ii) there is a $b \in B$, such that $b \notin ran(f)$; (iii) there is no function $g$ `greater than' $f$, that is, such that $ran(f) \subset ran(g)$ and $b \in ran(g)$.  

\end{definition}

\noindent
If one or more conditions (i)-(iii) do not hold, then one says that `$A$ is not smaller than $B$' (in symbols, $|A| \nless |B|$). Now we have the full statement of what we will call `Buzaglo's Principle':

\vspace{11pt}

\noindent
\textbf{Buzaglo's Principle (BuP)}. $|A|=|B| \leftrightarrow |A| \nless |B|$ and $|B| \nless |A|$.

\vspace{11pt}

\noindent
There are trivial cases of a set's not being smaller than another one: for instance, $\mathbb{R}$ is not smaller than $\mathbb{N}$, as, clearly, one can find an $f: \mathbb{R} \rightarrow \mathbb{N}$, such that there is no $n \notin ran(f)$. But the fact that $|\mathbb{N}| \nless |\mathbb{R}|$ is less trivial and helps to show that, by BuP: 

\[|\mathbb{N}|=|\mathbb{R}|\] 

\noindent
In fact, it can be shown that, for any uncountable set $A$: 

\[|\mathbb{N}|=|A|\] 

\noindent
so, BuP is one further example of a `one-cardinality' theory, which assigns the same cardinality to all \textit{infinite sets}. 

Although the motivations behind BuP are different from those behind Peano's early conception of infinite cardinalities, BuP shows that there are, in principle, ways to sanction the adoption of just one infinite cardinality, which may be as legitimate as those underlying the `bijection method' which sanctions the existence of Cantor's transfinite. Of course, whether, and in what sense, BuP (as much as Peano's `single-cardinality' conception, for that matter) is an interesting and mathematically fruitful principle is another question, which is besides the scope of this paper.

\pagebreak

\bibliography{Bib1}

\bibliographystyle{apalike}

\pagebreak

\tableofcontents

\end{flushleft}

\end{document}